\documentclass[12pt,a4paper]{amsart}
\usepackage[utf8]{inputenc}
\usepackage{amsmath}
\usepackage{amsfonts}
\usepackage{amssymb}

\usepackage{float}
\usepackage{subfig}

\usepackage{graphicx}
\usepackage{caption}

\usepackage{hyperref}

\newtheorem{thm}{Theorem}[section]

\newtheorem{lem}[thm]{Lemma}

\def\HH{\mathbb{H}}

\def\ZZ{\mathbb{Z}}
\def\CC{\mathbb{C}}
\def\RR{\mathbb{R}}
\def\QQ{\mathbb{Q}}
\def\NN{\mathbb{N}}

\def\tt{\Sigma_{1,1}}

\def\F2{\ZZ * \ZZ}
\def\aut{\text{Aut}(\F2)}
\def\out{\text{Out}(\F2)}
\def\gl2{\mathrm{GL}(2, \ZZ)}
\def\sl2{\mathrm{SL}(2, \ZZ)}
\def\slr{\mathrm{SL}(2, \RR)}
\def\slc{\mathrm{SL}(2, \CC)}

\def\hc{\tilde{\Sigma}}
\def\isom{\mathrm{isom}(\HH)}

\def\tr{\text{tr\,}}

\def\lab{\ell_{\alpha\beta^{-1}}}

\title{Rank two free groups and integer points on real  cubic surfaces}

 \author[McShane]{Greg McShane }
 \thanks{The author thanks ICERM, Brown University and the CNRS for support whilst preparing this article.}

\begin{document}

\begin{abstract} 
Counting integer  points on the Markoff cubic 
is closely related to questions in hyperbolic geometry.
In a previous work with Igor Rivin we investigated 
the regularity of the geodesic length function
for a punctured torus. 
Here we extend this work to the three holed sphere.

\end{abstract} 

\maketitle

\section{Introduction}

The rank two free group $\F2$ 
arises as the fundamental group of 
exactly four non-homeomorphic surfaces:
the three-holed sphere $\Sigma(0,3)$, 
the one-holed torus $\tt$, 
the two-holed projective plane (or one holed mobius band)  $C(0,2)$, 
and the one-holed Klein bottle $C(1,1)$. 
Of these the first two are orientable and the second two
are nonorientable. 
Furthermore $\tt$ enjoys the remarkable property
that every automorphism of its fundamental group
 is \textit{geometric} ie induced by a homeomorphism. 
Equivalently, every homotopy-equivalence $\tt \rightarrow \tt$
is homotopic to a homeomorphism. 
It is not surprising then  that 
much time and energy  has been 
devoted  to the geometric and topological 
properties of the  one holed torus. 
Here we will be corncerned with 
trying to extend the methods 
used for the one holed torus
to study the moduli spaces of the 
other surfaces.

  \begin{figure}[hb]
\centering
\includegraphics[scale=.4]{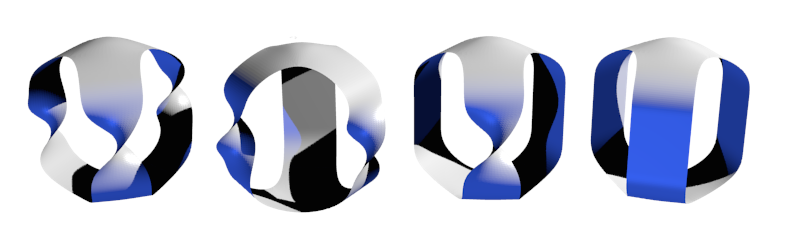} 
  \caption{From left to right- the one holed torus, one holed Klein bottle, one holed Mobius band, three holed sphere.}
  \label{pt graphic}
\end{figure}

\subsection{Hyperbolic surfaces}
Each of these four surfaces admits a (non unique) hyperbolic structure,
in fact there is a three dimensional deformation space of such structures.
Indeed it  has been known since the time of Fricke that,
as a consequence of Rieman's Uniformization Theorem, 
the three holed sphere and the one holed torus 
the moduli space of structures
can be identified with a semi algebraic subset of $\RR^3$. 
More precisely, 
one obtains a  hyperbolic structure  $X$ on $\Sigma$
from a discrete faithful representation $\rho$ of $\pi_1(\Sigma)$ 
 into $\isom$ such that $\HH/\rho(\F2)$  is homeomorphic to $\Sigma$.
 If $\Sigma$ is orientable one can lift $\rho$
 from a representation into  $\isom^+$,
  which is isomorphic to  $\text{P}\sl2$ to a representation, 
 $\tilde{\rho} : \F2 \rightarrow \text{SL}(2,\RR)$.
Doing this allows one to consider the trace  $\tr \tilde{\rho}(\gamma) \in \RR$ 
for any element $\gamma$ of the free group 
so that many problems can be reduced to questions of  linear algebra.
In particular one defines a character map 
$$\chi :  X \mapsto (\tr \tilde{\rho} (\alpha) , \tr \tilde{\rho} (\beta), 
\tr \tilde{\rho} (\alpha\beta)),$$
which provides an embedding of the moduli space of hyperbolic structures on $\Sigma$.
Goldman \cite{gold} 
studied the dynamical system defined by the 
action of the group of outer automorphisms 
of $\F2$ on the moduli space.
If $\phi$ is an automorphism then 
it acts on the representations 
$$ \tilde{\rho} \mapsto \tilde{\rho}  \circ \phi^{-1}$$
and evidently this induces an action 
on the points in the image of the embedding above.
It turns out that,
 by applying  Cayley-Hamilton theorem to $\sl2$,
it is easy to show that this action
is the restriction of polynomial diffeomorphisms  of $\RR^3$.
Later Goldman and his collaborators \cite{allofus} made  a similar study 
for the  two-holed projective plane  and the one-holed Klein bottle.
This latter work is quite delicate for two reasons.
Firstly, the  representation $\rho$  no longer lifts to $\slr$ 
but $\slc$
and secondly the automorphisms will in general  no longer be geometric
for these surfaces.

\subsection{Markoff triples}
An important aspect of the theory of the 
action of the group of outer automorphisms 
is that it admits an invariant function.
In the case of the pair of orientable surfaces this is 
$$ \kappa :  (x,y,z) \mapsto x^2 + y^2 + z^2 - x y z.$$
The level set $\kappa^{-1}(0)$ 
has five connected components - 
the singleton $(0,0,0)$
 and four codimension one manifolds homeomorphic to a $\RR^2$.
 The non trivial components can be identified with the 
 Teichmueller space of a once punctured torus.
 It is well known that there are  integer points in $\kappa^{-1}(0) \cap \RR_+^3$,
 for example $(3,3,3)$,
 and that they form a single orbit for the action of the outer automorphisms.
 From such a  point one obtains a \text{Markoff triple}, 
 that is a solution in positive integers  of the Markoff cubic 
 $$x^2 + y^2 + z^2 - 3x y z = 0 $$
 simply by dividing through by $3$.
 This map is \textit{natural} in that it is a conjugation
 of automorphisms of the level set $\kappa^{-1}(0)$
 with the automorphisms of the Markoff cubic.
So, 
since we will not be interested in arithemetic properties of these solutions,
we will also refer in what follows to elements of
 $\kappa^{-1}(0) \cap \NN^3 $ as Markoff numbers.
 An integer  which appears in a Markoff triple is called
 a \textit{Markoff number}.
 
 \subsubsection{Reduction theory}
 
 The fact that the Markoff triples form 
 a single orbit is a corollary of a classical result of Markoff 
 in reduction theory which we now explain briefly.
 We start by giving the Markoff triples an ordering,
 in the obvious way,
 using the sup norm on $\RR^3$ and proceed to show that
 any solution can be obtained from $(1,1,1)$
 by repeatedly applying generators of the  group of automorphisms 
 of the cubic.
This  automorphism group can be shown to be generated by:
 \begin{enumerate}
 \item sign change automorphisms $(x,y,z) \mapsto (x,-y,-z)$.
 \item coordinate permutations eg $(x,y,z) \mapsto (y,x,z)$.
 \item a Vieta flip $(x,y,z) \mapsto (x,y,3xy - z)$.
 \end{enumerate}
 Since Markoff triples are solutions in positive integers
 we will use just
 the \textit{Markoff morphisms} 
that is  the group generated 
 by permutations and the  Vieta flip.
 Given a Markoff triple 
 $(x,y,z) \neq (1,1,1)$  we may apply a permutation 
 so that  $x \leq y < z$
 which one can try to ``reduce" using the Vieta flip, 
 that is replacing the it  by the triple  $(x,y,3xy - z)$.
 We obtain a smaller solution provided $3xy < 2z$
 and this  inequality  holds 
 for every Markoff triple \textit{except}
 $(1,1,1)$.
 The reduction process gives rise
 to the structure of a rooted binary tree,
 the  \textit{Markoff tree},
 on the Markoff triples with $(1,1,1)$ as the \textit{root triple}.
 
  \begin{figure}[hb]
\centering
\includegraphics[scale=.3]{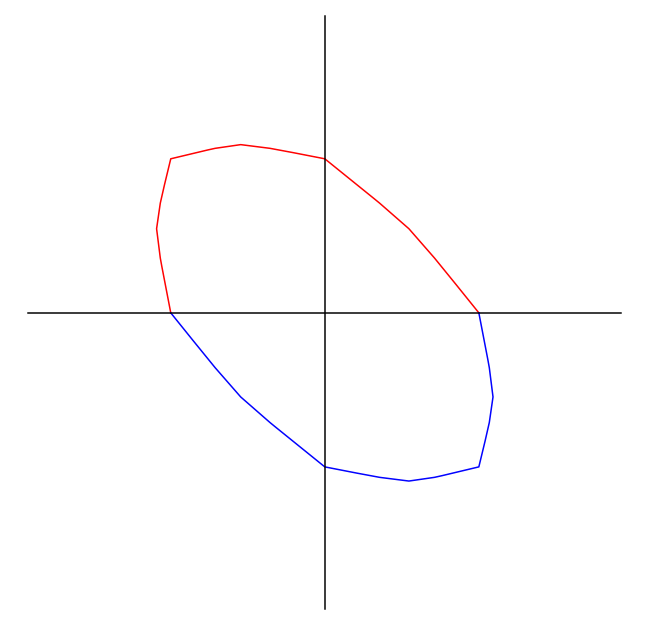} 
  \caption{Level set for the the length function $\ell$.}
\end{figure}
 
 \subsubsection{Counting integer points}
The question of counting for Markoff triples,
 that is for $R>0$ estimating the size of the set of Markoff triples
$(x,y,z)$ for which the sup norm is less than $R$, 
was first investigated by 
Gurwood  \cite{gu} in his thesis.
He  established an asymptotic formula using a correspondence
between Markoff and Farey trees.
Somewhat later Zagier \cite{za} obtained an improved  error term. 
At the time of writing
the  best result   is due to a Rivin and the author\cite{mcr}:
$$M(R) = C(\log R)^2 + O(\log R \log \log R),$$
for a numerical constant which has a geometric interpretation.

\begin{figure}[H]
\centering
\includegraphics[scale=.3]{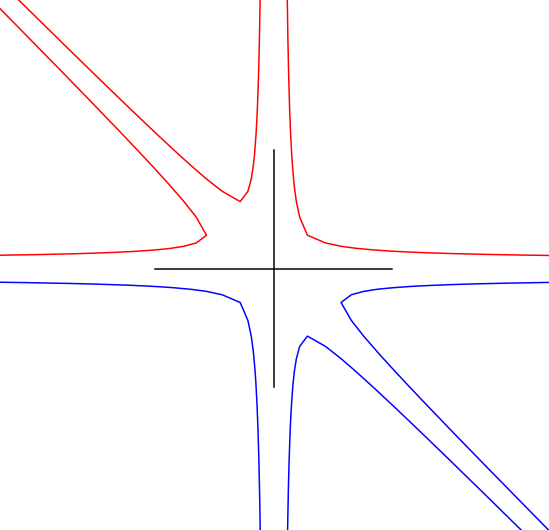} 
  \caption{Level set for the the length function $\ell$
  for the three punctured sphere.
  The scale is not the same as before as can be seen from the 
 representation of the black square cross representing the axes.}
  \label{three holed sphere graphic}
\end{figure}

The approach used in \cite{mcr} differs from that of Gurwood and Zagier
in that it uses ideas of H. Cohn  to interpret the question 
in terms of geometry of hyperbolic surfaces.
In \cite{cohn} Cohn 
 shows that the Markoff triple $(3,3,3)$ 
is the image under the character map  $\chi$ of the \textit{modular torus}
that is $X = \HH /\Gamma$ where $\Gamma$ is 
the commutator subgroup of $\text{PSL}(2, \ZZ)$.
It is easy to show that the modular surface is a 
6-fold cover of the modular surface $\HH /\text{PSL}(2, \ZZ)$.
Secondly, he shows that 
the permutations and the Vieta flips
are induced by automorphisms
and concludes that  the action of $\aut$
acts transitively on Markoff triples.
Now identifying $\F2$ with the fundamental group of the modular torus
one can show that 
 if $x$ is a  Markoff number then  there exists
a simple loop $\gamma$
such that $x = \tr \tilde{\rho}(\gamma)$.
It is well known that $\ell_X(\gamma)$,
the length of the closed geodesic  homotopic to $\gamma$ 
on the surface  $X = \HH/\rho(\F2)$,
 can be computed from the trace using the relation
 \begin{equation}\label{length from trace}
 2\cosh \ell_X(\gamma)/2 = |\tr \tilde{\rho}(\gamma)|,
  \end{equation}
so there is an explicit relation between Markoff numbers and
lengths of simple geodesics on the modular torus.
In  \cite{mcr} the simple closed curves are embedded 
as a subset of primitive elements in the integer lattice $\ZZ^2 \subset \RR^2$
and the length function $\gamma \mapsto \ell_X(\gamma)$ 
is shown, 
using a simple geometric argument,
 to extend to a proper  convex function $\ell: \RR^2 \mapsto \RR_+$.
 Thus the counting problem for Markoff numbers 
 can be settled by counting integer points in the set
 $$\ell^{-1}([0,R])  = R\ell^{-1}([0,1]) \subset \RR^2$$
 It is easy to visualize this set
 using a very short computer program
 see for example Figure \ref{three holed sphere graphic}.

 \subsection{Clebsch cubic}
 
 It is natural to study the distribution of integer points on 
 other affine surfaces $V_k(\RR)$ in the family, for $k \in \ZZ$
$$x^2 + y^2 + z^2 - 3x y z = k.$$
In one direction Ghosh,  Sarnak \cite{gs}
show that for almost all $k$ the Hasse Principle holds, 
namely that $V_k(\ZZ)$ is non-empty if
 $V_k(\ZZ_p)$ is non-empty for all primes $p$.
However they also show that there are infinitely
many k’s for which it fails. 
They give a proof  that the  Markoff morphisms again 
act on $V_k(\ZZ)$
 but with finitely many orbits rather than a single unique orbit.

There is a surface of particular interest in this family namely
$$x^2 + y^2 + z^2 - x y z = 20.$$
 Goldman has shown that this set can be identified with 
 the real points of the Clebsch cubic
 but what is of immediate interest to us is that
it  contains the integer point  $(2,2,-2)$.
This point is the image  under $\chi$ 
of a representation for which 
the corresponding hyperbolic surface $X$
 is a three punctured sphere.
Let us set $x=2$ so that 
\begin{equation}\label{parabolic line}
 16 = y^2 + z^2 - 2 y z   = (y -z)^2
 \end{equation}
ans we see that  there are infinitely many integer solutions on the lines $(2,t,t \pm 4), t \in \RR$ corresponding to an integer value of the parameter $t$.
An obvious consequence of this observation is 
that the growth rate of integer triples is quite different 
from that of the Markoff triples as it must be at least linear.

In fact these points
 fall into one of three orbits under the action of the Markoff morphisms
with the root of the corresponding Markoff tree at $(2,0,-4),(2,1,-3)$ and $(2,2,-2)$.
Using reduction theory Goldman \cite{gold} has shown that,
after possibly  applying a sign change automorphism,
any integer solution on this surface
is in the orbit of  one of these three triples under the Markoff morphisms.
Each of the  three solutions corresponds to hyperbolic structures on some surface
but only $(2,2-2)$ is non singular.
The triple $(2,1,-3)$  and $(2,0,-4)$ 
correspond to  structures on a punctured disc
each with  a single cone singularity 
of angle $2\pi/3$ and $\pi$ respectively.

If we choose one of the three root triples then,
using  the same process as was used to analyse 
the punctured torus in \cite{mcr},
we can vizualize the level sets of the corresponding length function.

\subsubsection{The three punctured sphere.}

The first triple we shall consider is $(2,2,-2)$ 
whose corresponding surface 
the three punctured sphere.

One sees in  Figure \ref{three holed sphere graphic}
that the set $\ell^{-1}([0,R]$ is quite different from 
that obtained in the case of the holed torus.
Though it is still symmetric,
that is invariant under half turn round the origin,
it is not compact for any $R>0$.
Moreover rather than being convex the function is ``piecewise" concave.
The six excursions to infinity correspond 
to oriented simple loops round one of the three punctures.
For each of these loops 
there is no minimiser in the homotopy class for  length,
informally we say that the geodesic has length zero,
and this is the reason for the infinite spike in this direction.

\subsubsection{The triple $(0,4,2)$.}

\begin{figure}[H]  
\centering
\includegraphics[scale=.6]{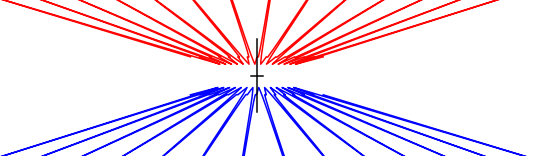} 
  \caption{Level set for the the length function $\ell$
 }\label{042graphic}
\end{figure}

The  triple $(0,4,2)$ 
corresponds to a  structure on a punctured disc
with  a single cone singularity  of angle $\pi$.
The boundary of the disc is a geodesic of strictly positive length and,
as before,
 we consider the puncture as being a geodesic of length $0$.
 The holonomy round the singularity is elliptic 
which  means that 
\begin{itemize}
\item 
the representation $\tilde{\rho}$ is not injective
\item 
 $(0,4,2)$ is fixed by an infinite cyclic subgroup of the Markoff morphisms.
 
\end{itemize}
 So we see  infinitely many  spikes as is illustrated in Figure \ref{042graphic}.

There is similar picture for $(0,3,1)$ which we will not reproduce here.

\subsubsection{The ``generic" case}

In general a hyperbolic structure on the three holed sphere 
will not have punctures or singularities
and the boundary holonomies will be hyperbolic.
In this case the minimum of the absolute value of 
the three numbers $x,y,z$ in the triple is strictly greater than $2$.
The level set for the length function resembles a six pointed star
with each of the prongs corresponding to 
a simple  loop  round one of the boundary components.

 \begin{figure}[H]
\centering
\includegraphics[scale=.4]{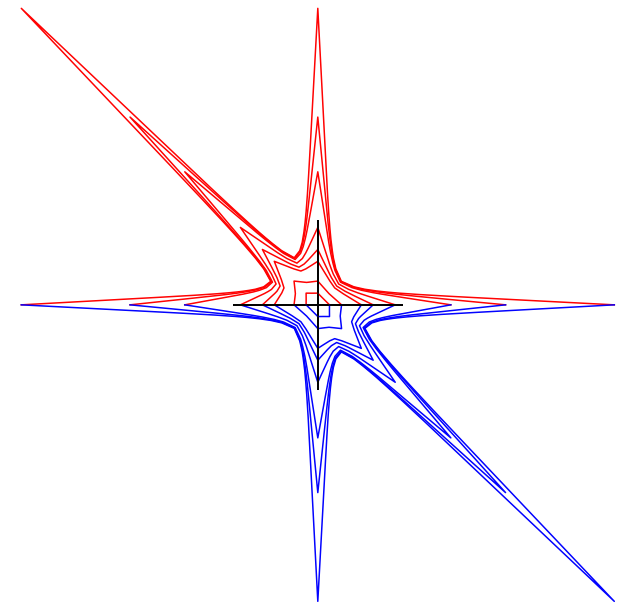} 
  \caption{Level sets of the length function $\ell$
  for  a collection of equilateral three holed spheres
   as  boundary length increases. 
 The level is compact, just as for the punctured torus,
 though it is no longer convex.}
  \label{longer and longer graphic}
\end{figure}

In Figure \ref{longer and longer graphic} we see
the level sets corresponding to 
several equilateral  three holed spheres.
Here by \textit{equilateral} we mean that the lengths of the
three boundary components are equal.
The lengths are actually $2.02,2.05,2.1,2.3,2.6,3,63$ and 6
and the figure illustrates  that as length increases 
the star shrinks 
and seams to be converging to a convex hexagon.

\section{Statement of results}

This work is concerned with the lengths of 
curves on oriented surfaces
we will treat the non orientable surfaces elsewhere.
The surfaces we treat here then are hyperbolic pairs of pants
which have three boundary curves of positive length.
The techniques apply to some extent 
to \textit{generalized pairs of pants}
where the boundary lengths can be zero 
or purely imaginary complex numbers,
which corresponds to a boundary being replaced by a cone singularity.
Our main result relates the shape of the level sets 
illustrated in the  figures of the previous paragraph
to geometric features of $X$ namely 
whether or not it has punctures and conical singularities.
What we are most  interested in is understanding 
why the level set of the length function 
is (locally) concave as opposed to being convex 
for a one holed torus.
Very roughly the difference comes down to the fact
that the minimiser for length in a homology class is 
a disconnected multi curve on the pair of pants 
whilst it is a single connected simple curve on the torus.

\begin{thm} \label{main theorem}
Let $X$ be 
a possibly incomplete hyperbolic structure on the three holed sphere.
Then the hyperbolic length function extends to 
a piecewise  concave function on $\RR^2$.

This extension vanishes
\begin{enumerate}
\item
at just the origin if $X$ has no conical singularity nor puncture
\item
on a  finite number of lines
corresponding to simple  loops round punctures
if $X$ has punctures but no conical singularity.
\item
a single line through the origin
if $X$ has a single  conical singularity 
but no cusps.
\item
countably many lines through the origin
if $X$ has a single  conical singularity
and at least one cusp.
\end{enumerate}

Further if $X$ has a cone singularity 
commensurable with $2\pi$ then 
the extension is invariant by a transvection.
\end{thm}

The points (1) and (2) 
follow from lower bounds for the length of a geodesic
which can be deduced from convexity.
The points (3) and (4)  follows from the fact that 
an  element of $\rho(\ZZ*\ZZ)$
representing the holonomy round the cone point
is elliptic.
Further if the angle at the  cone singularity is commensurable 
with $2\pi$ then this element is finite order $p$ say.
After choosing a loop round the cone singularity as the first
element of a  basis of the homology it is easy to show that
the level set will be invariant under:
$$ (x,y) \mapsto (x, y + px).$$

\subsubsection{Area and counting}

As in \cite{mcr} and \cite{Mir} it is possible 
to associate an area to $\ell^{-1}([0,1])$ 
in a natural way and we study how the  area varies with boundary length.
The two plots illustrate the main feature of this variation:
\begin{itemize}
\item  on the left we see  the variation of the level set 
corresponding to triples  $(2.3,2.3,t)$ as t varies from 4 to $2^{20}$
\item 
on the right we see a plot  for $(t,2.3,10)$ as t varies from 2.3 to 2.00243.
We see the spike growing, to accomodate this growth 
we show only half the level set and have changed the aspect ratio of the image.
\end{itemize}

 \begin{figure}[H]
  \centering
  \subfloat  
{\includegraphics[scale=.25]{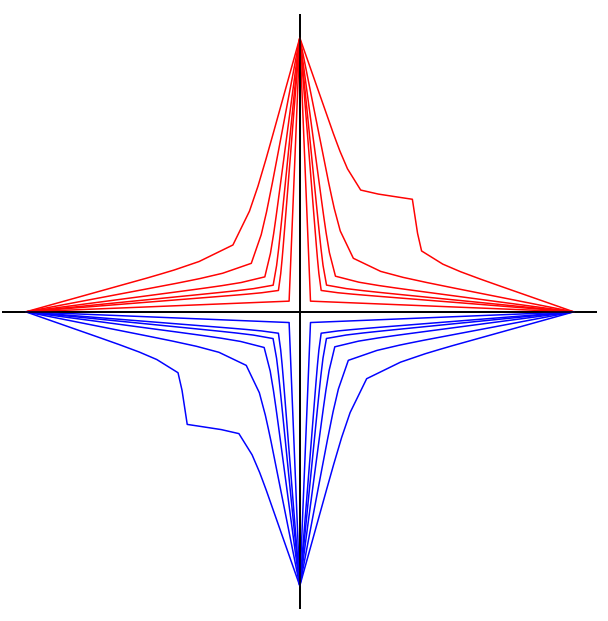} 
}
 \hfill
  \subfloat
  {\includegraphics[scale=.5]{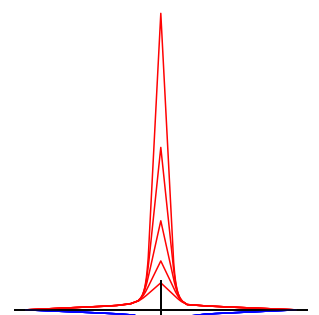} 
}
  \caption{Variation of the level sets as 
  one of the boundary length increases (on the left) 
  or decreases (on the right).
   }
\end{figure}

%
%
%

The fact that   $\ell^{-1}([0,1])$ 
may not be  compact  or even of finite area
means that counting integer points in 
the associated Markoff morphism orbit is very different from
the punctured torus \cite{mcr}.
Indeed whenever there is a parabolic element 
the growth of orbit will be essentially linear.
This is easy to see for the Clebsch cubic
from equation (\ref{parabolic line}).

\subsubsection{Non orientable surfaces}

A similar analysis can be carried out for the pair of non orientable surfaces  with fundamental group isomorphic to $\F2$. 
The level set of the length function  no longer bounds a convex region 
and appears to have a fractal nature.
\begin{figure}[H]
\centering
\includegraphics[scale=.3]{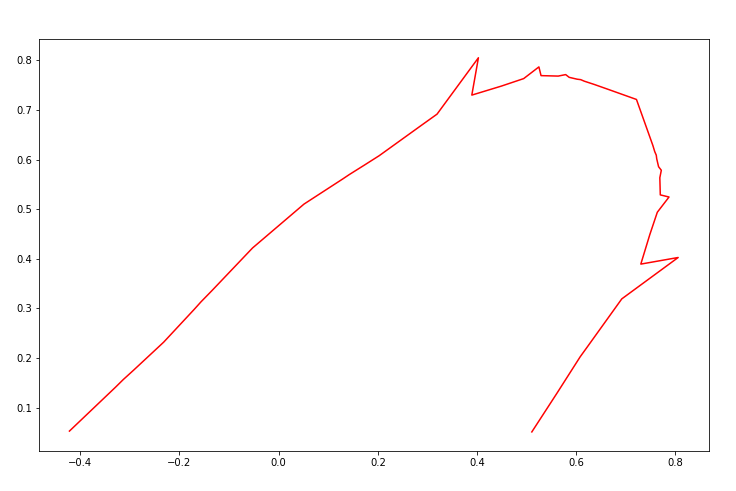} 
  \caption{Part of the level set of the length function of a non orientable surface.}
\end{figure}

\subsubsection{Markoff maps}

Following Bowditch, S.P. Tan and his co authors have studied
so-called \textit{Markoff maps} 
developping a theory analogous 
to Conway's topography for quadratic forms.
Some of our results can be interpreted in terms
of this theory. 
For example the excursions to infinity of the level sets 
correspond to the \textit{ends} of the Markoff map.
We will not however develop 
the relationship between our results
 and Markoff maps here .

\section{Automorphisms of $\F2$}
 
Let $\F2$ denote the free group on 2 generators
which we will denote $\alpha$ and $\beta$.
We think of  $\alpha$ and $\beta$ as simple loops 
meeting in a single point in a holed torus.
An element $\gamma \in \F2$ is called \textit{primitive}, 
if  there exists an automorphism $\phi$, such that 
$\gamma = \phi(\alpha)$. 
If $\phi(\delta) = \beta$, then $\gamma$ and $\delta$ are called
\textit{associated primitives}.

Following \cite{mcr}
let $\phi : \F2 \rightarrow \ZZ^2$
be the canonical abelianizing homomorphism.
The kernel of $\phi$ is a \textit{characteristic subgroup} 
that is it is $\aut$-invariant 
so there is a homorphism 
from $\aut$ to the automorphisms of $\ZZ^2$ 
namely $\gl2$.
In fact this homomorphism is surjective and 
the kernel is exactly the inner automorphisms
so that it induces an isomorphism between 
the group of outer automorphisms of $\ZZ*\ZZ$ and $\gl2$.

\begin{enumerate}
\item
If $\gamma$ is primitive then $\phi(\gamma)$ 
is a primitive element of $\ZZ^2$.
\item
$\aut$ acts transitively on primitive elements of $\F2$
and $\gl2$ acts transitively on primitive elements of $\ZZ^2$.

\end{enumerate}

\subsection{$\F2$ as a surface group}

The group $\F2$ is isomorphic to the fundamental group of exactly 4 surfaces, namely
\begin{itemize}
\item  the three holed sphere
\item  the one holed mobius band (the two holed projective plane)
\item  the one holed klein bottle
\item the one holed torus
\end{itemize}

\begin{figure}[H]
\centering
\includegraphics[scale=.5]{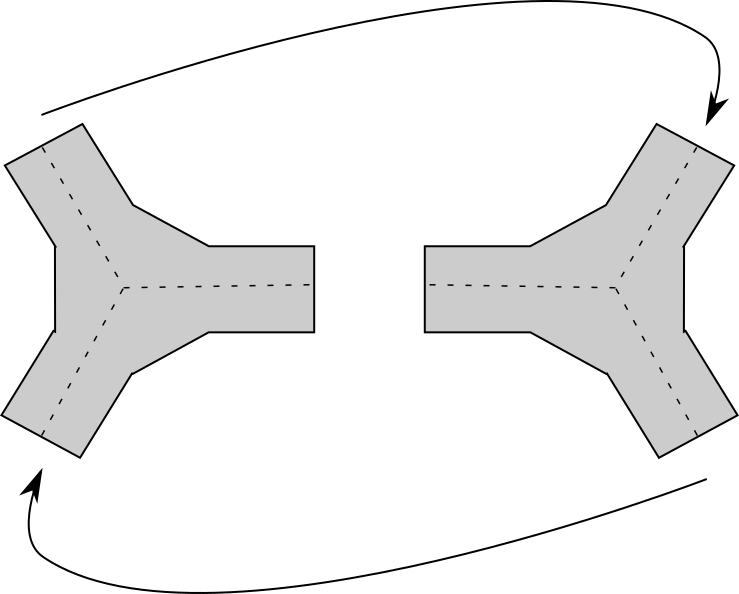} 
  \caption{A pair of tripods }
  \label{tripods}
\end{figure}

It is convenient to think of a topological surface with boundary as  a trivalent "fat graph".
A \textit{fat graph}  is a finite  connected
 graph equipped with a cyclic ordering on the half edges incident to each vertex.
To fat  graph, one can associate an oriented surface with boundary by replacing edges by thin oriented rectangles (ribbons) and vertices by disks, and pasting rectangles to disks according to the chosen cyclic orders at the vertices. 

In practice a fat graph is typically trivalent
and,
 if the fundamental group of the surface is $\F2$
 it  will have two vertices and three edges.
If we split each of these edges at its midpoint 
we have six half edges each adjacent to exactly one of the vertices.
A vertex with its three half edges is called a \textit{tripod}.
Each of the four surfaces in our list above
can be obtained as "twisted fat graphs"
from a pair of congruent  "fat tripods" as follows.
A fat tripod is a regular hexagon which has flaps attached to
three pairwise  non adjacent edges (see Figure ).
There are essentially four different ways to glue,
or more correctly sew, 
the tripods together 
depending on  whether the flaps are twisted once or not at all.
\begin{itemize}
\item if all the flaps are untwisted then one obtains a three holed sphere
\item if exactly one flap is twisted then one obtains a one holed mobius band
\item if exactly two flaps are twisted then one obtains a one holed klein bottle
\item if all the flaps are twisted then the result is a one holed torus.
\end{itemize}
Following Gilman and Keen \cite{GK enumerating}
we call the the three arcs along which the flaps are joined 
\textit{seams}.
 
Clearly, if $\Sigma$ is one of our four surfaces then 
one can obtain a convex hyperbolic structure on $X$ 
by replacing the tripods with a pair of congruent 
right angled hyperbolic hexagons.
When $\Sigma$ is a three holed sphere 
there are many good expositions of this procedure,
see for example \cite{Buser}, \cite{Beardon}.
In the case of  the punctured torus one obtains 
from  seams 
a triple of geodesic arcs perpendicular to the boundary.
These arcs are so-called \textit{ortho geodesics}
and this has structure has been quite extensively  studied 
by S.P. Tan and his co authors \cite{tan_torus}.
The hyperbolic structures on 
the one holed Mobius band and the one holed Klein bottle
are less well studied but appear in \cite{norbury} and \cite{allofus}.

\subsection{Topological realizations of automorphisms of $\F2$}

It is often useful to identify $\F2$ with the fundamental group of  the holed torus
and  think of $\out$, the outer automorphism group of $\F2$,
as being  the mapping class group of the holed torus.
The outer automorphism group is generated 
by the so-called \textit{Nielsen transformations}, which either permute the basis 
$\alpha, \beta$, or transform it into $\alpha \beta, \beta$, or
$\alpha \beta^{-1}, \beta$. 
Both these  latter transformations can be realized  topologically as  Dehn twists on the holed torus.

In what follows we will identify $\ZZ^2$ with the homology of the surface $\Sigma$
and  will call the cover $\hc \rightarrow \Sigma$ corresponding 
to $\ker \phi$ the \textit{homology cover or maximal abelian of $\Sigma$}.
 \begin{figure}[H]
\centering
\includegraphics[scale=.25]{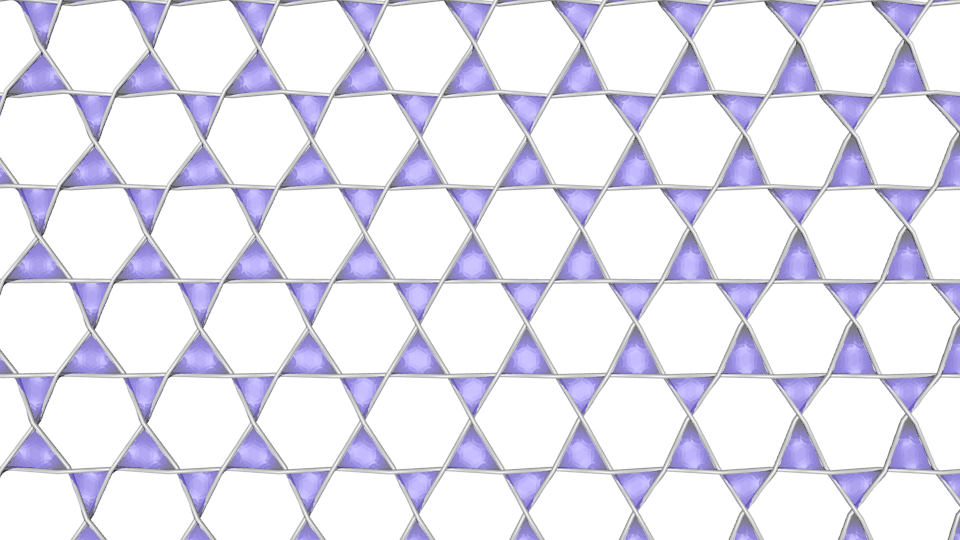} 
  \caption{Maximal abelian cover for the three holed sphere.
  The boundary loops lift to the grey curves,
  the patches appear triangular 
  but in fact are hexagons joined along sides 
  which are perpendicular to the plane.}
\end{figure}

There are many fine expositions of the action of $\out$
on the maximal abelian cover of the one holed torus.
Unfortunately the action does not lift well to 
the maximal abelian cover for the three holed sphere.
However, we will use a ``common spine" to transfer information.


\subsection{Spines and retractions}

Recall that a loop $\gamma$  is \textit{primitive} 
if there is an automorphism $\phi$ of the free group 
such that $\phi^{-1}(\gamma)$ is a generator (boundary loop).
It is well known that
 if $\alpha$ is a primitive loop  on a punctured torus 
 then  the unique geodesic freely homotopic  to $\alpha$ passes through  exactly two of the three Weierstrass points.
Recall that a \textit{Weierstrass point}
is a fixed point of the elliptic involution.
Our reasoning depends on transferring 
other well-known results for primitive geodesics
on the punctured torus to 
primitive geodesics a three holed sphere (pants).
This is possible because the punctured torus 
and the three holed sphere have spines
which are isomorphic -
we think of them as having a \textit{common spine}.

A \textit{spine} for a surface  is a 1-complex which 
is a strong deformation retract,
the dotted lines in Figure (\ref{tripods})
will glue together to make a spine for the surface.
Epstein and Bowditch used hyperbolic 
geometry to construct spines of surfaces 
with boundary or cusps.
 The spine for  a pair of pants
 consists of all the points $x$
in the complement 
such that  the distance to the boundary
is realised by at least two geodesic arcs.
To construct a spine for  the punctured torus
it is a little more complicated:
one first chooses a cusp region $H$
then the spine  consists of all the points $x$
in the complement of $H$
such that  the distance $d(x,H)$
is realised by at least two geodesic arcs.
It is easy to see that this does not depend on
the choice of $H$.
It is easy to see that 
the spine for a pair of pants 
is always a trivalent graph with exactly two vertices.
For a generic punctured torus this also the case.

We will use a spine for a punctured torus which 
is the union of two simple  loops which generate 
the fundamental group.
For the pair of pants our spine will be a 
\textit{figure eight curve}.

\begin{figure}[H]
\centering
\includegraphics[scale=.4]{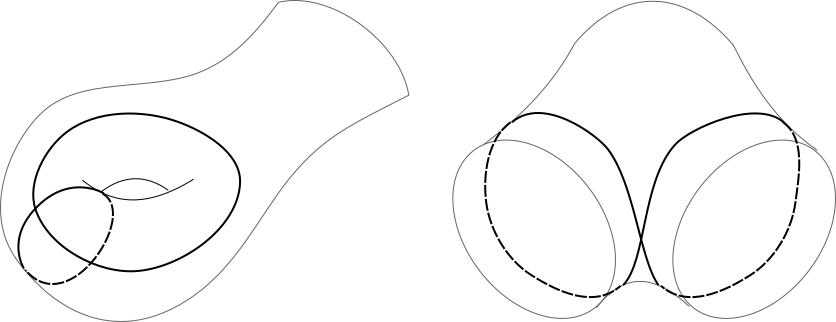} 
  \caption{Spines for surfaces. A pair of generating loops 
  in a one holed torus and a figure eight curve in a pair of pants.}
\end{figure}

\section{Non vanishing and lower bounds for length}

We begin by studying the case when there are no punctures.
The pants has three boundary components of
lengths $\ell_\alpha, \ell_\beta, \lab$.
Given a curve $\gamma$ such that 
$$\gamma  = m  \alpha + n \beta,$$
with $m,n$  non negative integers
then it is possible to do surgery 
at self intersection points of $\gamma$
to obtain a simple multi curve 
each of whose components is homotopic to
either  $\alpha$ or $\beta$.
This procedure reduces length and so 
$$\ell_\gamma  \geq  m  \ell_\alpha + n \ell_\beta,$$

Suppose that the surface has a single puncture 
so that, with the notation above,  $\ell_\alpha = 0$.
We think of the surface as being obtained 
from an essential  sub surface of the thrice punctured sphere.
That is $X \subset Y:= \HH/\Gamma(2)$ 
(where as usual $\Gamma(2)$ is the principal congruence subgroup)
so that $X$ inherits its complex structure from the inclusion 
and a hyperbolic structure via the Poincar\'e metric.
By the Ahlfors-Pick-Schwarz Lemma
the length of $\gamma \subset X$
is bounded below by the length of $\gamma \subset Y$
so it suffices to bound the length on the thrice punctured sphere.

\begin{figure}[H]
\centering
\includegraphics[scale=.3]{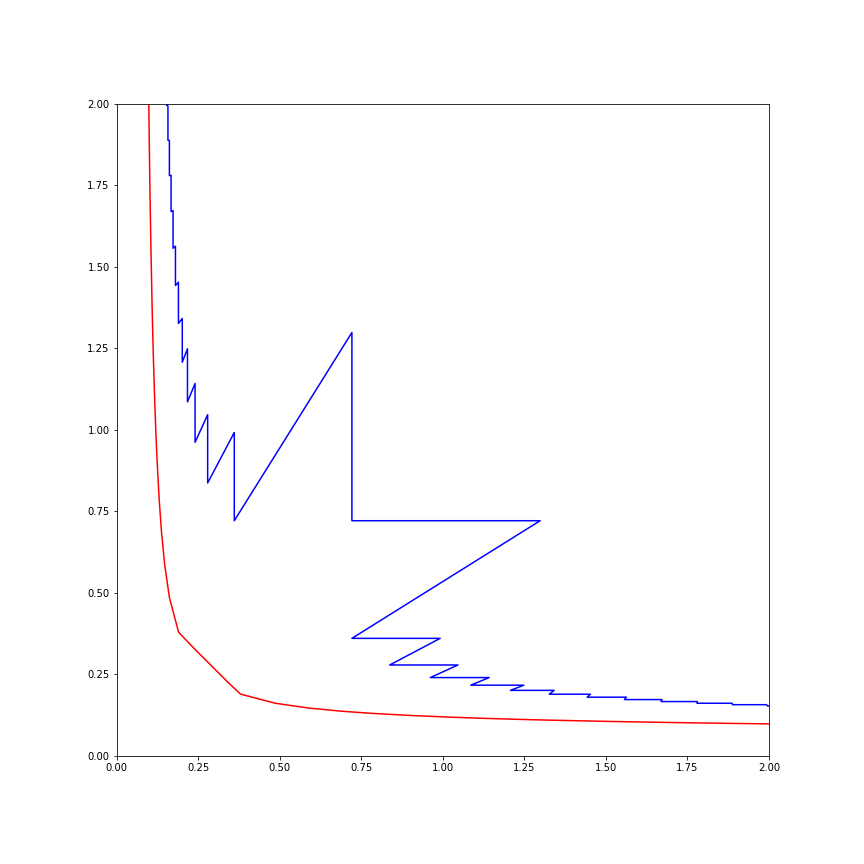} 
  \caption{Hyperbolic length versus our estimate.}
  \label{compare lengths}
\end{figure}

The length of  $\gamma$ is obtained from the trace 
of a primitive element of the fundamental group which we denote 
$C$.
We choose generators for the fundamental group
which are a pair of loops round $\alpha$ and $\beta$
denoted  $A$ and $B$ respectively.
Gilman \cite{JG} has made a study of the word in $A,B$ that represents $C$.
Relabelling boundaries 
 if necessary we may assume $m>n$ 
 and then this word takes the form
$$C = A^{m_1}B\ldots A^{m_r}B$$
where the $m_i \in \{ [m/n], [m/n] + 1\}$.
One can then argue that the length of $\gamma$
is approximately
$$n \times \ell_{\gamma_1}$$
where $\gamma_1$ is the loop representing the homotopy class of 
$ A^{m_1}B$.
The length of the $\gamma_1$ can be calulated from the trace of
the matrix
$$\begin{pmatrix}
1 & 2 \\ 0 & 1
\end{pmatrix}^{m_1}
\begin{pmatrix}
1 & 0  \\ 2 & 1
\end{pmatrix}
= 
\begin{pmatrix}
1 +  4m_1&  2m_1  \\ 2 & 1
\end{pmatrix} \in \Gamma(2).
$$
so that 
$
\cosh \ell_X(\gamma_1)/2 = 2m_1 + 1,
$
and 
$$\ell_X(\gamma_1) \geq 2n\log( 2m_1).$$
It is not difficult to see that 
if $m > n$ and $m/n \leq K$,
so that $n \geq m/K$,
then 
$$2n\log \left(2\left[\frac{m}{n} \right] \right)
 \geq   \frac{\log(2)}{K} m .$$
We have plotted the graph of this estimate 
and the graph of the catual length in Figure \ref{compare lengths}.
This is enough to show the non vanishing of the extension
as in (2) of Theorem \ref{main theorem}.

\section{Geodesic length and convexity for the holed torus}


\subsection{Convexity for the holed torus}

Let $\gamma \in \F2$ be some closed curve 
and $\ell_X(\gamma)$ 
the minimum over the lengths  wrt the hyperbolic 
structure $X$ of all closed 
curves  in the homotopy 
class of $\gamma$.
In fact the  minimum is attained and is 
lthe length of the unique (oriented)  closed geodesic
in the homotopy  class.

\begin{thm}
The  length function is convex on $H_1(X,\QQ)$
and so extends continuously (as a convex function) to $H_1(X,\RR)$.
\end{thm}

Since the proof is short and instructive we shall reproduce it now.

Let $X$ be a holed torus equipped with a hyperbolic
structure. 
Then, the shortest multicurve representing a non-trivial homology
class $h$ is a  simple closed geodesic  if $h$ is a primitive homology class 
(that is, not a multiple of another class), and a multiply covered geodesic otherwise. 
Further, the shortest multicurve representing h is unique.

Begin by defining a valuation $\ell$ on the first homology  with integral coefficients, 
where $\ell(h)$ is just the length of the shortest multicurve representing h. 
The valuation of the trivial homology class is defined to be 0. 
By the previous paragraph $\ell$ statisfies
$$ \ell(nh) = n\ell(h), \forall n \in \NN.$$
Moreover if $h,g$ are not commensurable 
then $\ell$  satisfies the strict triangle inequality 
\begin{equation} \label{triangle inequality}
\ell(h + g) < \ell(h) + \ell(g).
\end{equation}

This is because the union of the shortest
multicurves representing  h and g is not embedded,
that is there is a transverse intersection somewhere.
 Hence, the shortest multicurve representing h + g 
 is strictly shorter than the union.

\section{Concavity for the three holed sphere}

To prove concavity we have to prove that,
in a certain sense,
the opposite of the triangle holds
that is:
\begin{equation} \label{anti triangle inequality}
\ell(h + g) > \ell(h) + \ell(g).
\end{equation}
Whilst the triangle inequality followed as a 
direct consequence of the fact that 
a multiple of a simple geodesic
in the class $h + g$ was  the unique minimiser,
the inequality (\ref{anti triangle inequality})
is be more difficult to prove.
Firstly, it is not true for any choice of $h$ and $g$
and this is why we can only prove concavity piecewise.
Secondly, 
although our proof is also  based on eliminating transverse 
intersections to reduce length the argument is more delicate.

So we begin by showing that (\ref{anti triangle inequality})
is true for a few very simple configurations.
We will then give an induction argument 
 to show that it is true in general.

\subsection{Eliminating intersections}
The key point in proving (\ref{triangle inequality}) 
is the fact that the multi geodesics 
representing $h$ and $g$
had a transverse intersection somewhere 
that could be eliminated, or \textit{surgered out},
to find a shorter representative for $h+g$.
This is iterated and must terminate when there 
are no longer any transverse intersections.

Given $h$ and $g$  there are multi geodesic representatives
on  the three holed sphere 
but, unless they are boundary loops,  they are no longer simple.
Typically they have self intersections on the seams of the three holed sphere that Gilman and Keen \cite{GK cutting}\cite{GK enumerating}
 call \textit{essential intersections}.
In a certain sense the opposite of 
what happened in the holed torus occurs.
Beginning with  $k\gamma$ the multi geodesic
representing $h+g$ and we eliminate an essential intersection 
to get a pair of multi geodesics which 
have (fewer transverse intersections)
smaller total length than $\ell(k\gamma) = k \ell(\gamma)$.
Iterating this process we should be able to show the required inequality.
For example the figure eight curve has a single self intersection
and after surgery it gives rise to a pair of curves
homotopic  to distinct boundary curves (which are disjoint).
In the next section we treat this example more formally.

\subsection{The figure eight curve.}

Let us begin by considering the triple 
$\alpha,\beta$ and $\alpha\beta$
where, as usual
$\alpha,\beta, \alpha\beta^{-1}$  
are simple loops round the boundary components
and $\alpha\beta$ is a figure of eight curve.
The figure of eight curve has a single essential intersection point
and there are two loops meeting at this point
one homotopic to $\alpha$ and the other $\beta$
( in fact the surface retracts onto this figure eight curve
and we will need this later.)
Clearly since neither of these loops is smooth they 
are longer than the closed geodesic  in the same homotopy class
and since length of $\alpha\beta^{-1}$ is the sum of these loops
we have:
\begin{equation} \label{anti triangle inequality1}
\ell(\alpha\beta) > \ell(\alpha) + \ell(\beta).
\end{equation}

Let's take this further consider the pair of loops
$$(\alpha\beta)^n\beta, (\alpha\beta)^{n}\alpha$$
so that in the abelianisation one has
\begin{equation*}
(n+1) \alpha\beta  = (\alpha\beta)^n\beta + (\alpha\beta)^{n}\alpha.
\end{equation*}
The  loops $(\alpha\beta)^n\beta$ and $(\alpha\beta)^{n}$
retract onto the figure eight curve  
$\gamma$
and by choosing our retraction carefully
we may suppose that this map is of degree $n+1$.
Now the  loops on the $\gamma$
are not smooth so they are stricly longer than the geodesics 
in the  the same  homotopy class  and so
\begin{equation} \label{anti triangle inequality2}
(n+1) \ell(\alpha\beta) 
> \ell((\alpha\beta)^n\beta) + \ell( (\alpha\beta)^{n}\alpha).
\end{equation}
We define the 
\textit{positive cone determined  by a pair of oriented loops}
to be the subset of the homology consisting 
of linear  combinations with positive coefficients.
More generally if we have
\begin{equation*} \label{xx}
k \alpha\beta  = m  \gamma_1 + n \gamma_2
\end{equation*}
with $\gamma_1,\gamma_2$
primitive loops in the positive cone  generated by 
the oriented loops $\alpha$ and $\beta$ in the homology
then,
since the surface retracts onto 
the figure eight graph determined by the  $\alpha\beta$
each of the loops $\gamma_i$ is homotopic to 
a piecewise geodesic curve on $\gamma$.
So, essentially by the same argument,
the reverse triangle inequality is still true, that is 
\begin{equation} 
k\ell(\alpha\beta) > m\ell(\gamma_1) + n\ell(\gamma_2),
\end{equation}
where $k$ is the degree of the retraction 
from the multi curve $m\gamma_1 \cup n\gamma_2$.

\subsection{Convexity in general}

To prove convexity we decompose 
our loop $\gamma$ into a pair of generator loops
 $\alpha',\beta'$,
 which are not necessarily boundary loops,
but  so that $\gamma$ is still a product
 \begin{eqnarray*}
 \gamma = \alpha'\beta'.
 \end{eqnarray*}
We do this by  splitting 
at an essential self intersection then writing $\gamma_i$ 
as words in these.
With such a decomposition the reverse triangle inequality follows
since $\alpha$ and $\beta$ minimise length in their respective
homotopy classes 
by applying the argument of the previous paragraph.

More formally :

\begin{thm}\label{decomposition}
Let $\gamma_1,\gamma_2$ be a pair of primitive geodesic loops
such that there are positive integers $m,n \geq 0$  with
$$k \gamma  = m  \gamma_1 + n \gamma_2.$$
in the abelianisation.
Then there is a positive integer $r$ such that 
$rk\gamma$  can be decomposed into
$rm$ loops homotopic to $\gamma_1$
and $rn$ loops homotopic to $\gamma_2$
\end{thm}

\subsubsection{Very local concavity}

We can prove  a version of Theorem \ref{decomposition}:
decomposition for pairs $\gamma_1, \gamma_2$ 
in a (non uniformly)  small neighborhood  of  the loop $\gamma$
and with $c=1$.
From this we can deduce concavity on the neighborhood
and by a local to global argument on the whole positive quadrant.
By a \textit{small neighborhood   of $\gamma$} we mean
we mean those primitive $\gamma'$ such that 
 $\gamma/\|\gamma\|_2$ and  $\gamma'/\| \gamma'\|_2$
are close (where, as usual $\|.\|_2$ denotes the Euclidean norm).
In fact, if $\gamma = p \alpha + q \beta$ with $p,q \geq 0$ coprime integers
then this neighborhood will be the positive cone determined by the 
loops
\begin{eqnarray*}
a \alpha + c \beta \\
b \alpha + d \beta
\end{eqnarray*}
where $a/c$ and $b/d$ are Farey neighbors of each other and  $p/q$.

Our argument uses the notion of a
\textit{staircase curve},
that is a curve
which is contained in the grid
$ (\RR \times (\ZZ + \frac12) )
 \cup  ((\ZZ + \frac12) \times \RR) \subset \RR^2$
 note that this grid   is wholly contained
  in the punctured plane $\RR^2 \setminus \ZZ^2$.

\begin{figure}
\centering
\includegraphics[scale=.3]{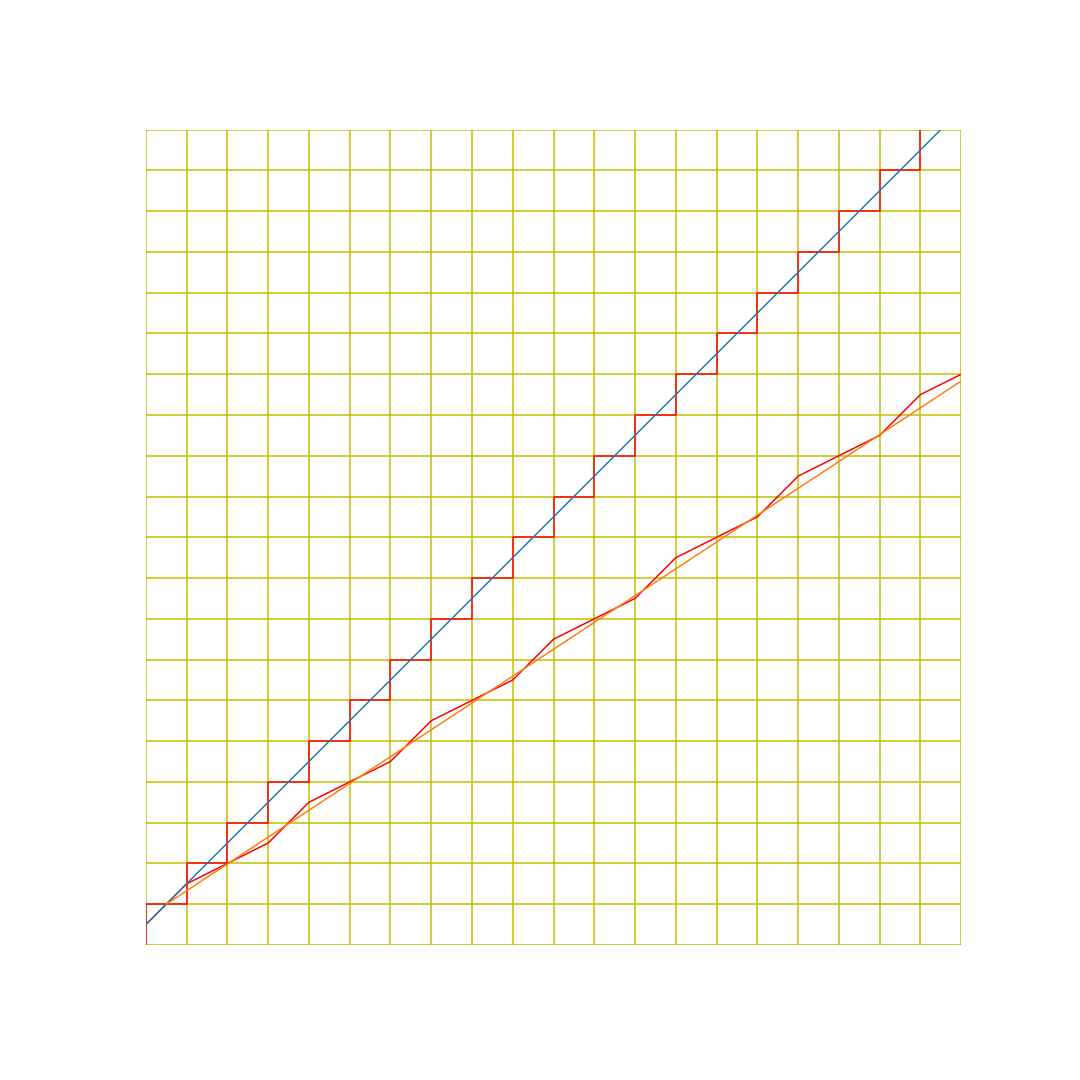} 
  \caption{Transfering concavity}
  \label{transfer}
\end{figure}

The punctured plane $\RR^2 \setminus \ZZ^2$
can be identified with the maximal abelian cover of the punctured torus:
if  $\Gamma\simeq \ZZ^2$ denotes the group of deck transformations
then the punctured torus is obtained as a quotient:
$$\left(\RR^2 \setminus \ZZ^2 \right)/ \Gamma.$$
There is a $\Gamma$-equivariant retraction 
$H_t$ from the punctured plane to a grid of 
horizontal and vertical lines (see Figure (\ref{transfer})).
So this grid is in fact the set of lifts of 
 a spine for the once punctured torus.

Recall  that a  \textit{Weierstrass point} on the punctured torus 
is a fixed point of the \textit{elliptic involution}
that is the involution of 
$\left(\RR^2 \setminus \ZZ^2 \right)/ \Gamma$
induced by $\RR^2 \rightarrow \RR^2,\, v \mapsto -v$.
It is easy to see that the set of lifts of these points are respectively
$\ZZ^2 + (\frac12,0)$ , $\ZZ^2 + (0,\frac12)$ and
$\ZZ^2 + (\frac12,\frac12)$.
So the grid contains lifts of the three Weierstrass points - 
these  occur at the vertices of the grid
and the midpoints of the horizontal and vertical edges.
We may also  identify the grid with the lift of a spine 
for the pants. 
We quite simply choose a figure eight curve and a retraction onto 
it such that each of the three essential intersections
is fixed.
Thus we think of the grid as being a ``common spine"
for the maximal abelian cover of both the punctured torus
and the pair of pants.
This allows us to visualise
essential intersections of geodesics
as  lifts of Weierstrass points on the grid.

Consider the line
$$\hat{\gamma}
= \left\{  \begin{pmatrix}
t \\ t + \frac 12
\end{pmatrix} , t \in \RR
\right\}$$
this curve is important as its
retraction to the spine of $H_1(\hat{\gamma})$
is a staircase curve which 
coincides with a lift of the figure eight curve.
The figure eight curve crosses 
the three seams of the pants
but only one gives rise to an essential self intersection.
The corresponding Weierstrass point 
is the one at the vertex of the grid.
Doing surgery at each vertex 
yields the horizontal and vertical lines
which are lifts of curves homotopic to 
boundary components of the pants.
Any other straight line 
$\hat{\zeta}_1$
avoiding the  integer lattice
and with positive gradient
retracts onto a ``staircase curve" 
$H_1(\hat{\zeta}_1)$ on 
the grid  and so decomposes 
as a sequence of horizontal  and vertical edges
from which we obtain the partition as 
in the previous subsection.


Now  consider $p ,q>0$ coprime integers. 
It is well known that 
$p/q$ is the descendant of two \textit{Farey neighbors}
that is a pair $a/c,b/d \in \QQ_+$
with $ad -bc= 1$,
such that 
$$ \frac{p}{q} = \frac{a+b}{c+d}.$$
This means that we can think of 
$\begin{pmatrix}
p \\  q
\end{pmatrix}$
as being  the image of
$\begin{pmatrix}
1 \\  1
\end{pmatrix}$
under the matrix
$$M 
= \begin{pmatrix}
a & b\\ c & d
\end{pmatrix}  
 \in \sl2.$$
So under the linear transformation 
 $f_M:v \mapsto Mv$ our line
 $\hat{\gamma}$ is mapped to $f_M(\hat{\gamma})$
 which is  a line of gradient $p/q$
avoiding the half integer lattice.
Note that the image of the staircase curve $H_1(\hat{\gamma})$
is not a staircase curve but it is a union of line segments
and homotopic to $f_M(\hat{\gamma})$.
Now for any $p'/q'$ such that 
$$ \frac{a}{c}  \leq \frac{p'}{q'}	 \leq \frac{b}{d}$$
the image under $f_M^{-1}$ of the ray through 
$\begin{pmatrix}
p' \\  q'
\end{pmatrix}$
is a line  $\hat{\zeta}$ with positive gradient
which  retracts  onto
the staircase curve  $H_1(\hat{\zeta})$
as in the previous paragraph.
Thus $\zeta$ decomposes as a sequence 
lifts of sub loops of the figure eight curve.
Finally  $f_M$  takes 
$\hat{\zeta}$ to a line of gradient $p'/q'$
and $H_1(\hat{\zeta})$ 
to a curve which is a union of line segments
between lifts of Weierstrass points/essential intersections
and so they  yield the decomposition of the loop $f_M(\zeta)$
as subloops of the loop determined by
$f_M(\gamma)$ on the pair of pants.

This procedure is illustrated for $p=3,q=2$
in Figure (\ref{transfer}).

Finally there is some accounting to be done 
if $p'/q'$ and $p''/q''$ are in the small neighborhood
bounded by the Farey ancestors of $p/q$ 
then there exists $m,n$ such that 
$$k\begin{pmatrix}
p \\q
\end{pmatrix}
= m\begin{pmatrix}
p' \\q'
\end{pmatrix}
+ n\begin{pmatrix}
p'' \\q''
\end{pmatrix}
$$
where $k$ is the determinant of
$$\begin{pmatrix}
p' & p''\\ q' & q''
\end{pmatrix},$$
and this proves the theorem at least  locally.

\subsubsection{Bigger neighborhoods} 
Now we can expand the convexity to a bigger 
neighborhood  by 
inheriting the neighborhoods of the  Farey ancestors 
$a/c$ and $b/d$ of $p/q$.
By the preceding paragraph there is a neighborhood
of $a/c$ for which each curve retracts onto 
a sequence of sub loops of the corresponding 
closed geodesic.
Again by the preceding paragraph each of these 
sub loops contracts onto a sub loop of the 
closed geodesic associated to $p/q$.

Iterating this process, 
and keeping track of the  ``accounting"
one can prove the theorem on the whole of the positive quadrant.

\end{document}